\documentclass[11pt]{article}
\usepackage[a4paper]{anysize}\marginsize{3.0cm}{3.0cm}{1.5cm}{1.5cm}
\pdfpagewidth=\paperwidth \pdfpageheight=\paperheight
\usepackage{graphicx}
\usepackage[colorlinks=true, linkcolor=blue, anchorcolor=blue, citecolor=blue, filecolor=blue, menucolor= blue, urlcolor=blue]{hyperref}
\usepackage{pgf,tikz-cd}
\usetikzlibrary{arrows}
\usepackage{amsfonts,amssymb,amsthm,amsmath,euscript}
\usepackage{ltablex}
\usepackage{caption, booktabs}
\usepackage{subcaption}
\usepackage{graphicx,color}
\usetikzlibrary{cd}
\usepackage[T1]{fontenc}
\usepackage{breqn}

\theoremstyle{plain}
\newtheorem{thm}{Theorem}[section]
\newtheorem{theorem}[thm]{Theorem}

\newtheorem{lemma}[thm]{Lemma}
\newtheorem{proposition}[thm]{Proposition}
\newtheorem{corollary}[thm]{Corollary}

\theoremstyle{definition}
\newtheorem{definition}[thm]{Definition}
\newtheorem{remark}[thm]{Remark}
\newtheorem{example}[thm]{Example}

\newtheorem{thevarthm}[thm]{\varthmname}

\newenvironment{varthm*}[1]{\trivlist\item[]{\bf #1.}\it}{\endtrivlist}

\newcommand\cali{\mathcal I}
\newcommand\calo{\mathcal O}

\newcommand\calc{\mathcal C}
\newcommand\Z{\mathbb Z}

\newcommand\C{\mathbb C}
\newcommand\E{\mathbb E}

\renewcommand\P{\mathbb P}
\newcommand{\eps}{\varepsilon}
\newcommand\wtilde\widetilde
\newcommand\newop[2]{\def#1{\mathop{\rm #2}\nolimits}}
\newop\codim{codim}
\newop\rank{rank}

\makeatletter
\newcommand{\subjclass}[2]{%
  \let\@oldtitle\@title%
  \gdef\@title{\@oldtitle\footnotetext{#1 \emph{Mathematics subject classification.} #2}}%
}
\newcommand{\keywords}[1]{%
  \let\@@oldtitle\@title%
  \gdef\@title{\@@oldtitle\footnotetext{\emph{Key words and phrases.} #1.}}%
}
\makeatother

\title{Companion varieties for root systems\\ and Fermat arrangements}
\author{Roberta Di Gennaro, Giovanna Ilardi,  Rosa Maria Mir\'o-Roig 
\\ Tomasz Szemberg, Justyna Szpond}
\date{\today}
\subjclass{2010}{13A15, 13C70, 14C20, 14E05, 14J70, 14N20}

\begin{document}
\definecolor{uuuuuu}{rgb}{0.26,0.26,0.26}
\maketitle
\thispagestyle{empty}

\begin{abstract} Unexpected hypersurfaces are a brand name for some special linear systems. They were introduced around 2017 and are a field of intensive study since then. They attracted a lot of attention because of their close connections to various other areas of mathematics including vector bundles, arrangements of hyperplanes, geometry of projective varieties. Our research is motivated by  what is now known as the BMSS duality, which is a new way of deriving projective varieties out of already constructed ones. The last author coined the concept of companion surfaces in the setting of unexpected curves admitted by the $B_3$ root system. Here we extend this construction in various directions. We revisit the configurations of points associated to either root systems or to Fermat arrangements and we study the  geometry of the associated varieties. In the case of configurations of points associated to  root systems, the geometry of their companions is also described.
\end{abstract}

\section{Introduction}
\label{sec: intro}
   In the present note we study companion varieties of unexpected hypersurfaces associated to some symmetric configurations of points in projective spaces. Companion varieties were introduced by the last author in \cite[Section 4.4]{Szp19} in the course of her study of Togliatti-type surfaces.
   
   It is a classical problem in algebraic geometry to determine the dimension of linear systems. A lot of attention was given to linear systems with imposed base loci, i.e., of the form
   \begin{equation}
   \label{eq: h0}
     H^0(X;L\otimes \cali(Z)),       
   \end{equation}
   where $L$ is a positive (e.g. ample or very ample) line bundle on a smooth variety $X$ and $Z$ is a subscheme of $X$. Even in the very simple situation when $X$ is the projective plane, $L$ is the line bundle $\calo_{\P^2}(d)$ for $d>0$ and $Z$ is a zero-dimensional subscheme of $\P^2$ the problem of computing the dimension of the vector space in \eqref{eq: h0} is far from being completely understood. Indeed, the two open conjectures: one due to Nagata (1959) \cite{Nag59} and the other the SHGH-Conjecture package due to Segre (1969), Harbourne (1986), Gimigliano (1987) and Hirschowitz (1989) (see \cite{DKMS16shgh} for introduction and recent progress) provide sufficient evidence for this claim. In any case, it is well-understood and clear that a single point, or a fat point scheme concentrated in a single general point impose independent conditions on homogeneous polynomials of any fixed degree in a projective space of arbitrary dimension. With this in mind, it came as a surprise when Cook II, Harbourne, Migliore and Nagel \cite{CHMN}, inspired by earlier results due to Di Gennaro, Ilardi and Vall\`es \cite{DIV14}, announced that a single general fat point might impose fewer conditions than expected on the linear system of homogeneous polynomials with assigned base loci. The elements arising this way are called \emph{unexpected hypersurfaces}, see Section \ref{sec: unexpected} for precise definitions. 
   
   The unexpected hypersurfaces are interesting not only because they violate the naive conditions count when determining the dimension of a linear system. In a private communication to the last author Igor Dolgachev suggested that, under additional positivity conditions, there should be a relation between the existence of unexpected hypersurfaces and varieties with defective osculating spaces. This has been indeed established in \cite[Section 4.2]{Szp19c} for a configuration of points determined by the $B_3$ root system. The linear system of quadrics vanishing at configuration points determines, after passing to the blow up of $\P^2$ in these points, a morphism to $\P^5$ whose image is a surface $S$ such that at every point $P$ of $S$, there is a hyperplane in $\P^5$
   tangent to $P$ to order $2$ (in other words: cutting out on $S$ a curve which passes through $P$ with multiplicity at least $3$). It has been observed additionally, that there is another surface $S'$, which we call a companion surface of $S$, which exhibits also interesting geometrical properties. This example motivates our present work.
   
   The general yoga, in the basic version taken from \cite{Szp19c}, is as follows. Assume that there is a set of points $Z$ in $\P^N$ which admits a unique unexpected hypersurface $H_{Z,P}$ of degree $d$ and multiplicity $m$ at a general point $P=(a_0:\cdots:a_N)\in \P^N$. Let 
   $$F_Z((a_0:\cdots:a_N),(x_0:\cdots:x_N))=0$$
   be a homogeneous polynomial equation of $H_{Z,P}$. Let $g_0,\ldots,g_T$ be a basis of the vector space
   $[I(Z)]_d$ of homogeneous polynomials of degree $d$ vanishing at all points of $Z$. Under some mild hypotheses the unexpected hypersurface $H_{Z,P}$ comes from a bi-homogeneous polynomial $F_Z((a_0:\cdots:a_N),(x_0:\cdots:x_N))$ of bidegree 
   $(m,d)$, see \cite[Section 4]{HMNT} and \cite[Proposition 1]{Matrixwise}.  Indeed, $F_Z$ can be written in a unique way as a combination
   \begin{equation}\label{eq: F as combination of generators}
      F_Z=h_0(a_0:\cdots:a_N)g_0(x_0:\cdots:x_N)+\cdots+h_T(a_0:\cdots:a_N)g_T(x_0:\cdots:x_N)
   \end{equation}
   where  $g_0(x_0:\cdots:x_N), \ldots ,g_T(x_0:\cdots:x_N)$ are homogeneous polynomials of degree $d$ forming the basis mentioned above and $h_0(a_0:\cdots:a_N),  
    \ldots ,h_T(a_0:\cdots:a_N)$ are homogeneous polynomials of degree $m$.
   Therefore, there are two rational maps naturally associated to equation \eqref{eq: F as combination of generators}
   $$\varphi: \P^N\ni (x_0:\cdots:x_N)\dashrightarrow (g_0(x_0:\cdots:x_N):\cdots:g_T(x_0:\cdots:x_N))\in\P^T$$
   and 
   $$\psi: \P^N\ni (a_0:\cdots:a_N)\dashrightarrow (h_0(a_0:\cdots:a_N):\cdots:h_T(a_0:\cdots:a_N))\in\P^T.$$
   The image of $\psi$ is the \emph{companion variety} of the image of $\varphi$. The purpose of this note is to study properties of companion varieties and relations between them. Turning to details in Sections \ref{sec: root systems} and \ref{sec: Fermats} we work in a more general situation where $g_0,\ldots,g_T$ form a basis of a vector subspace in $[I(Z)]_d$ large enough to write down $F_Z$ in the form in \eqref{eq: F as combination of generators}. It might also happen that the coefficients $h_0,\ldots,h_T$ are linearly dependent. In this case we work rather with their linearly independent subset in order to avoid dealing with degenerate subvarieties. Thus it is possible (and it actually happens) that $\varphi$ and $\psi$ are mappings to projective spaces of different dimensions.

   We focus here on three natural generalizations of the aforementioned $B_3$ configuration. In what follows $H$ is always the pull-back of the hyperplane bundle under the appropriate blow up.
   
   First, results from \cite{HMNT} suggest that all
$B_n$ configurations give rise to unexpected hypersurfaces in $\P^{n-1}$. In case of $n=4$ there is also a unique unexpected hypersurface and we determine its explicit equation in \eqref{eq: unexpected B4} and in slightly different generators in \eqref{eq: unexpected cone for F4}. It turns out that the unexpected surface for $B_4$ is also unexpected for $F_4$,
   a non-crystallographic root system. The other difference when compared to the $B_3$ root system is that now the degree of the unexpected surface equals its multiplicity in the general point, i.e., the unexpected surfaces are cones.
\begin{varthm*}{Theorem $B_4$/$F_4$}
   Let $Y$ be the blow up of $\P^3$ at points in the $F_4$ root system with exceptional divisors $E_1,\ldots,E_{24}$. Then the linear system $4H-\sum E_i$ embeds $Y$ into $\P^{11}$ as a smooth threefold $X$ of degree $40$. The companion threefold $X'$ is in this case isomorphic to $X$.
\end{varthm*}   
   This is proved in Proposition \ref{prop: threefold of deg 40} and Remark \ref{rem: B4 and F4 self companion}.
   
   Secondly, we study the $H_3$ root system. Similarly as $B_3$ for $m=3$, it admits an unexpected curve of degree $m+1$ with multiplicity $m$ at the general point, this time with $m=5$.
\begin{varthm*}{Theorem $H_3$}
   Let $Y$ be the blow up of $P^2$ at points in the $H_3$ root system with exceptional divisors $E_1,\ldots,E_{15}$. Then the linear system $6H-\sum E_i$ embeds $Y$ as a smooth aCM surface $X$ of degree $21$. The companion surface $X'$ is not linearly normal embedded into $\P^{11}$.
\end{varthm*}   
   This result is proved in Proposition \ref{prop: surface of segree 21} and Remarks \ref{rem: H3 companion} and \ref{rem: H3 bicompanion}.
   A very surprising feature of this example is that the companion surface $X''$ of $X'$ is not $Y$ again but it is rather its projection.
   
   Finally, we study companion surfaces related to Fermat arrangements. The $B_3$ configuration belongs also to this family. Our main results in this direction are Theorem \ref{thm: Fermat Z(Fm,3)}, Proposition \ref{prop_image_Fm} and Remark \ref{rem: 5.8}. The precise statements are slightly too technical to quote them here.
  
   We work over the field of complex numbers $\C$.

\section{Unexpected hypersurfaces}
\label{sec: unexpected}
   The ground-breaking work \cite{CHMN} by Cook II, Harbourne, Migliore and Nagel introduced the concept of unexpected curves. This notion was generalized to arbitrary hypersurfaces in the subsequent article \cite{HMNT} by Harbourne, Migliore, Nagel and Teitler.
\begin{definition}\label{def:unexpected hypersurface}
   We say that a reduced set of points $Z\subset\P^N$ \emph{admits an unexpected hypersurface} of degree $d$
   if there exists a sequence of non-negative integers $m_1,\ldots,m_s$ such that for general points $P_1,\ldots,P_s$
   the zero-dimensional subscheme $P = m_1P_1+\cdots +m_sP_s$ fails to impose independent conditions on forms
   of degree $d$ vanishing along $Z$ and the set of such forms is non-empty. In other words, we have
   $$h^0(\P^N;\calo_{\P^N}(d)\otimes I(Z)\otimes I(P))>
     \max\left\{0, h^0(\P^N;\calo_{\P^N}(d)\otimes I(Z))-\sum_{i=1}^s\binom{N+m_{i}-1}{N}\right\}.$$
\end{definition}

Let us illustrate the above definition with an example.

\begin{example}\label{ex1} We consider the extended Fermat arrangement of planes in $\P^3$: 
$$F=xyzw(x^4-y^4)(x^4-z^4)(x^4-w^4)(y^4-z^4)(y^4-w^4)(z^4-w^4)$$
and we denote by $Z\subset \P^3$ the set of 28 points dual to the linear factors of $F$. We have:
$$\begin{aligned}
   I(Z)=&(xyz,xyw,xzw,yzw,\\
   &xy(x^4-y^4),xz(x^4-z^4),xw(x^4-w^4),yz(y^4-z^4),yw(y^4-w^4),zw(z^4-w^4)).
\end{aligned}$$   
Using {\em Macaulay2} \cite{M2}, we easily check that there is one and only one unexpected surface of degree 6 containing $Z$ and three general points $P_1$, $P_2$, $P_3$ of multiplicities $m_1=m_2=m_3=4$.
\end{example}

   Following \cite[Definition 2.5]{ChiantiniMigliore19} we introduce also the following notion.
\begin{definition}[Unexpected cone property]
   Let $Z$ be a finite set of points in $\P^N$ and let $d$ be a positive integer. We say that $Z$ has the \emph{unexpected cone property} $\calc(d)$, if for a general point $P\in\P^N$, there exists an unexpected (in the sense of Definition \ref{def:unexpected hypersurface}) hypersurface $S_P$ of degree $d$ and multiplicity $d$ at $P$ passing through all points in $Z$. 
\end{definition}

\begin{example}
We consider the set $Z\subset \P^3$ of 28 points introduced in Example \ref{ex1}. Using again {\em Macaulay2} \cite{M2},  we check that   $Z$ has the unexpected cone property $\calc(6)$ since for a general point $P\in\P^3$, there exists an unexpected surface $S_P\subset \P^3$ of degree $6$ and multiplicity $6$ at $P$ passing through all 28 points in $Z$.
\end{example}
\section{The BMSS duality}
\label{sec: duality}

For the sake of completeness let us explain what we understand for BMSS duality a notion which comes from 
\cite{BMSS} and \cite[Section 4]{HMNT}.

We consider integers $0<m\le d\in \Z$, a set $Z\subset \P^N$ of points and a general point $P=(a_0:a_1:\cdots :a_N)\in \P^N$. We assume that there is a unique hypersurface $H_Z\subset \P^N$ of degree $d$ containing $Z$ and having multiplicity $m$ at $P$.
The hypersurface $H_Z$ is defined by a bihomogeneous form
$$F_Z(\underline{y},\underline{x}):=F_Z((y_0:y_1:\cdots :y_N),(x_0:x_1:\cdots :x_N))\in \C[\underline{y},\underline{x}].$$
of bidegree $(t,d)$ with $t\ge m$ (see \cite[Lemma 3.1(d)]{HMNT}).

Note that the form  must be bihomogeneous since if it were not, changing the projective coordinates of each of the two set of variables, the form would vary. Moreover, under mild hypotheses on $Z$ the bidegree of $F_Z(\underline{y},\underline{x})$ is $(m,d)$ i.e. $d=m$ (see, for instance,   \cite[Proposition 1]{Matrixwise} for the general case and \cite[Section 4]{HMNT} for special cases where much more can be said).
It is important to point out that  the bidegree of $F_Z(\underline{y},\underline{x})$ is not always $(m,d)$. In fact, in \cite[Theorem 6]{BMSS} the authors considered a set $W\subset \P^3$ of 31 points  with an expected quartic surface $S_W\subset \P^3$  having a general point $Q\in \P^3$ of multiplicity  3 and defined by a bihomogeneous polynomial $F_W(\underline{y},\underline{x})$ of bidegree $(5,4)$.

The BMSS duality allows us to view the bihomogeneous polynomial $F_Z(\underline{y},\underline{x})\in \C [\underline{y}][\underline{x}]$ of bidegree $(d,m)$ as a polynomial of degree $m$ in the indeterminates $\underline{x}$ with coefficients the point $P\in \P^N$ of multiplicity $m$, i.e., a family of hypersurfaces in the variables $\underline{x }$ parameterized by $P$ (the unexpected hypersurface):
$$F_{Z,P}(\underline{x}):=F_Z((a_0:\cdots:a_N),\underline{x})\in \C[a_0,\ldots ,a_N][\underline{x}]$$ 
or, as a homogeneous polynomial of degree $d\ge m$ in the variables $\underline{y}$ with coefficients a point $Q=(b_0:\cdots :b_N)\in \P^N$, i.e. a family of hypersurfaces in the variables $\underline{y}$ parameterized by $Q$:
$$F_{Z,Q}(\underline{y}):=F_Z(\underline{y},(b_0:\cdots:b_N))\in \C[b_0,\ldots ,b_N][\underline{y}].$$

The BMSS duality establishes when the following assertions are true: 
\begin{itemize}
    \item  The tangent cone of $F_{Z,P}(\underline{x})$ at $Q$ coincides with  the tangent cone of $F_{Z,Q}(\underline{y})$ at $P$.
    \item  $F_{Z,P}(\underline{x})$ has multiplicity $m$ at a general point $Q$ and  $F_{Z,Q}(\underline{y})$ has also multiplicity $m$ at a general point $P$.
\end{itemize}

In next section we discuss in detail an example which comes from the free arrangement of planes in $\P^3$ determined by linear factors of:
$$xyz(x+y)(x-y)(x+z)(x-z)(x+w)(x-w)(y+z)(y-z)(y+w)(y-w)(z+w)(z-w)$$
as well as an example which comes from the root system $H_3$.

\section{Companion varieties for root systems}
\label{sec: root systems}
   In \cite[Section 3]{HMNT} the authors study unexpected hypersurfaces related to root systems. These are finite sets of vectors in a vector spaces satisfying certain number of conditions, see \cite{Humphreys} for introduction, motivation and basic properties of this important notion. For us root systems are just a source of interesting configurations of points in projective spaces and we don't use explicitly any of many properties they enjoy. A vector $v$ in the root system in the vector space $V$ determines a point in the projective space $\P(V)$, which corresponds to the line spanned by $v$. As the set $Z$ one takes all points in $\P(V)$ determined in this way by all vectors in the root system.
   
   The study in \cite{HMNT} was motivated by \cite{DIV14}, where the first example of an unexpected curve coming from the $B_3$ root system was found. Computer experiments carried out in \cite{HMNT} provide the following list of root systems admitting a unique unexpected hypersurface. In Table \ref{tab:root systems with 1 unexpected} the number $N$ stands for the dimension of the ambient projective space $\P^N$ the number $d$ for the degree of an unexpected hypersurface and $m$ for the multiplicity in the general point.
\begin{table}[ht]
    \centering
    \begin{tabular}{c|c|c|c}
       root system & N & d & m\\
       \hline
       $B_3$ & 2 & 4 & 3\\
       \hline
       $B_4$ & 3 & 4 & 4\\
       \hline
       $D_4$ & 3 & 3 & 3\\
       \hline
       $F_4$ & 3 & 4 & 4\\
       \hline
       $H_3$ & 2 & 6 & 5\\
       \hline
       $H_4$ & 3 & 6 & 6
    \end{tabular}
    \caption{Root systems determining a unique unexpected hypersurface}
    \label{tab:root systems with 1 unexpected}
\end{table}

   As \cite[Section 3]{HMNT} and \cite[Section 2]{BMSS} study the $B_3$ case in detail, we turn  our  attention to the $B_4$ root system and to the $H_3$ root system. We will encounter a couple of new phenomena in the course:
   \begin{itemize}
       \item in the case of the $B_4$ root system, after passing to the appropriate subspace of $[I(B_4)]_4$, we obtain a variety with an isomorphic companion variety;
       \item in the case of the $H_3$ root system we obtain the companion surface which is not aCM and which is embedded in a space of lower dimension than the initial variety.
   \end{itemize}
\subsection{The \texorpdfstring{$B_4$}{B4} root system}
The set $Z(B_4)$ consists of $16$ points, which can be assigned the following coordinates:
$$\begin{array}{llll}
P_{1} = [1:0:0:0] & P_{2} = [0:1:0:0] & P_{3} = [0:0:1:0] & P_{4} = [0:0:0:1]
\end{array}$$
$$\begin{array}{lll}
P_{5} = [1:1:0:0] & P_{6} = [1:0:1:0] & P_{7} = [1:0:0:1] \\
\end{array}$$
$$\begin{array}{lll}
P_{8} = [0:1:1:0] & P_{9} = [0:1:0:1] & P_{10} = [0:0:1:1] \\
\end{array}$$
$$\begin{array}{lll}
P_{11} = [1:-1:0:0] & P_{12} = [1:0:-1:0] & P_{13} = [1:0:0:-1] \\
\end{array}$$
$$\begin{array}{lll}
P_{14} = [0:1:-1:0] & P_{15} = [0:1:0:-1] & P_{16} = [0:0:1:-1] 
\end{array}$$
   The saturated ideal $I(B_4)$ is generated by
$$xyz, xyw, xzw, yzw,$$   
and
\begin{equation}\label{eq: generators of B4 in deg 4}
xy(x^2-y^2), xz(x^2-z^2), xw(x^2-w^2), yz(y^2-z^2), yw(y^2-w^2), zw(z^2-w^2).  
\end{equation}

A basis for the vector space $[I(B_4)]_4$ is given by including
with the degree $4$ elements listed in \eqref{eq: generators of B4 in deg 4}, the following elements:
$$xyzw, x^2yz, xy^2z, xyz^2, x^2yw, xy^2w, xyw^2, x^2zw, xz^2w, xzw^2, y^2zw, yz^2w, yzw^2.$$   
   The unexpected cone,
   written down in terms of this basis has then the following equation
\begin{align}\label{eq: unexpected B4}
\begin{split}
   F=&3ad(c^2-b^2)x^2yz+ 3bd(a^2-c^2)xy^2z+ 3cd(b^2-a^2)xyz^2+\\  
   &3ac(b^2-d^2)x^2yw+ 3bc(d^2-a^2)xy^2w+ 3dc(a^2-b^2)xyw^2+\\
   &3ab(d^2-c^2)x^2zw+ 3cb(a^2-d^2)xz^2w+ 3db(c^2-a^2)xzw^2+\\
   &3ba(c^2-d^2)y^2zw+ 3ca(d^2-b^2)yz^2w+ 3da(b^2-c^2)yzw^2+\\
   &cd(d^2-c^2)xy(x^2-y^2)+ bd(b^2-d^2)xz(x^2-z^2)+ bc(c^2-b^2)xw(x^2-w^2)+\\
   &ad(d^2-a^2)yz(y^2-z^2)+ ac(a^2-c^2)yw(y^2-w^2)+
   ab(b^2-a^2)zw(z^2-w^2).
\end{split}
\end{align}
   It is immediately clear, that the generator $xyzw$ is not involved in the equation. The second property we observe is also quite surprising.
\begin{proposition}\label{prop:F4 as base locus}
   All unexpected cones for the $B_4$ root system, have in their base locus $8$ additional points:
$$\begin{array}{lll}
Q_1=[1:1:1:1] & Q_2=[1:1:1:-1] & Q_3=[1:1:-1:1] \\ Q_4=[1:1:-1:-1] &
Q_5=[1:-1:1:1] & Q_6=[1:-1:1:-1] \\  Q_7=[1:-1:-1:1] & Q_8=[1:-1:-1:-1].&
\end{array}$$
\end{proposition}
\proof  
   The vanishing of $F$ at the given points is easy to check. That the base locus is not larger can be verified by computing equations of unexpected cones for some specific values of $a,b,c,d$.
\endproof
   The points in $B_4$ together with the points $Q_1,\ldots,Q_8$ give the root system $F_4$.
\begin{corollary}
   The unexpected cones for $B_4$ and $F_4$ are the same.
\end{corollary}
   This observation has interesting implications for the companion threefolds. First, let us note, that 
   a basis for $[I(F_4)]_4$ is given by the following forms:
\begin{equation}\label{eq: generators for F4}
\begin{array}{lll}
   m_0=xy(x^2-y^2), & m_1=xy(z^2-w^2), & m_2=xz(x^2-z^2), \\
   m_3=xz(y^2-w^2), & m_4=xw(x^2-w^2), & m_5=xw(z^2-y^2), \\
   m_6=yz(y^2-z^2), & m_7=yz(x^2-w^2), & m_8=yw(y^2-w^2), \\
   m_9=yw(x^2-z^2), & m_{10}=zw(z^2-w^2), & m_{11}=zw(x^2-y^2).
\end{array}
\end{equation}
   Using these generators the equation \eqref{eq: unexpected B4} becomes
\begin{align}\label{eq: unexpected cone for F4}
\begin{split}
   F=&
   cd(d^2-c^2)xy(x^2-y^2)+
   3cd(b^2-a^2)xy(z^2-w^2)+
   bd(b^2-d^2)xz(x^2-z^2)+ \\
   &
   3bd(a^2-c^2)xz(y^2-w^2)+
   bc(c^2-b^2)xw(x^2-w^2)+ 
   3bc(a^2-d^2)xw(z^2-y^2)+ \\
   &
   ad(d^2-a^2)yz(y^2-z^2)+
   3ad(c^2-b^2)yz(x^2-w^2)+ 
   ac(a^2-c^2)yw(y^2-w^2)+ \\
   &
   3ac(b^2-d^2)yw(x^2-z^2)+
   ab(b^2-a^2)zw(z^2-w^2)+
   3ab(d^2-c^2)zw(x^2-y^2).
    \end{split}
\end{align}   
   Let $\varphi:\P^3\dashrightarrow\P^{11}$ be the rational map defined by the generators in \eqref{eq: generators for F4}.
\begin{proposition}\label{prop: threefold of deg 40}
   The image of $\varphi$ is a smooth threefold of degree $40$.
\end{proposition}   
\proof
   Let $\sigma:X\to\P^3$ be the simultaneous blow up of each of $24$ points in $F_4$ with the exceptional divisor $\E$ (which splits into $24$ projective planes, one over each of the points blown up) and as usual let $H=\sigma^*\calo_{\P^3}(1)$. 
    Since $\{ p\in \P^3 \mid m_{i}(p)=0 \text{ for } i=0,\ldots , 11 \}=\{ P_1,\ldots ,P_{16},Q_1,\ldots ,Q_8\}$ and since the forms $m_i$ are unions of planes which at each base point $p$ have no common tangent, the linear system  $L=4H-\E$ is base point free and it defines a morphism onto its image 
    $\varphi_{L}:X\to\P^{11}$ which lifts the map
    $$\varphi:\P^3\dashrightarrow\P^{11} \quad \quad (x:y:z:w)\mapsto (m_0:m_1:\cdots :m_{11}).$$
 Therefore, we have the commutative diagram
 $$ \begin{tikzcd} 
 Y \arrow[d,"\sigma "] \arrow[dr,"\varphi _{L}"] &\\
 \P^3 \arrow[r ,dashrightarrow, "\varphi"] & \P^{11}.
 \end{tikzcd}
 $$
 Let us call $X$ the image  of $\varphi$.
 To check that the map  $\varphi : \P^3 \dashrightarrow Y$ has degree 1 and even more it is a birational map onto its image, we use the packages "Cremona"   (\cite{CremonaSource}, \cite{CremonaArticle}) and "RationalMaps" \cite{RationalMaps} of {\em Macaulay2} \cite{M2}. In fact, we have the following code block:
 
 \vskip 4mm
\noindent\texttt{loadPackage "Cremona"}\\
\texttt{loadPackage "RationalMaps"}\\
\texttt{R = QQ[x,y,z,w]}\\
\texttt{I =  ideal(x*y*(x$\wedge$2-y$\wedge$ 2),x*y*(z$\wedge$2-w$\wedge$2),x*z*(x$\wedge$2-z$\wedge$2),x*z*(y$\wedge$2-w$\wedge$2),}\\
\texttt{ x*w*(x$\wedge$2-w$\wedge$2),x*w*(y$\wedge$2-z$\wedge$2),y*z*(y$\wedge$2-z$\wedge$2),y*z*(x$\wedge$2-w$\wedge$2),}\\
\texttt{ y*w*(y$\wedge$2-w$\wedge$2),y*w*(x$\wedge$2-z$\wedge$2),z*w*(z$\wedge$2-w$\wedge$2),z*w*(x$\wedge$2-y$\wedge$2))}\\
\texttt{degree rationalMap gens I}\\    
(Using the package "Cremona" we get that the degree is 1)\\
\texttt{S=QQ[a\_0..a\_((numgens I)-1)]}\\
\texttt{L = apply(numgens I, j -> (gens I)\_{(0,j)})}\\
\texttt{varphi = map(R, S, L)}\\
\texttt{isBirationalOntoImage(varphi)}\\
(Using the package "RationalMaps" we get that the map is birational onto its image.)
   \vskip 4mm
   Once we know that $\varphi$ is a birational map onto its image we compute the degree of $Y$:  $degree (X)=L^3$. The smoothness of $X$ follows from the Jacobian criterion and a straightforward computation which shows that outside $\{ P_1,\ldots ,P_{16},Q_1,\ldots ,Q_8\}$ the rank of $\begin{pmatrix}  \frac{\partial m_{i}}{\partial x} & \frac{\partial m_{i}}{\partial y} &\frac{\partial m_{i}}{\partial z} & \frac{\partial m_{i}}{\partial w} \end{pmatrix} _{i=0,\ldots ,11 }$ is 4.
   
   It is worthwhile to point out that the Jacobian dual criterion as stated in \cite[Theorem 2.18 and 3.2]{DHS} or \cite[Theorem 1.4]{RS} gives us a computer free proof of the birationality of $\varphi $. In fact, the ideal $I=(m_0,m_1,\ldots,m_{11})\subset R:=\C [x,y,z,w]$ has a linear presentation:
   $$
   \cdots \longrightarrow  R(-5)^{16} \xrightarrow{\psi} R(-4)^{12} \longrightarrow R \longrightarrow R/I \longrightarrow 0
   $$
   with 
 
    \begin{equation*}
   \psi=\left(\begin{array}{cccccccccccccccc}
   
  -z & 0 & 0 & 0 & 0 & 0 & -w & 0 & 0 & 0 & 0 & 0 & 0 & 0 & 0 & 0 \\

   y & -y & 0 & 0 & 0 & 0 & 0 & 0 & 0 & 0 & -w & 0 & 0 & 0 & 0 & 0 \\
   
   0 & x & -y & 0 & -z & 0 & 0 & 0 & 0 & 0 & 0 & 0 & -w & 0 & 0 & 0 \\
   
   -y & 0 & x & -y & 0 & -z & 0 & 0 & 0 & 0 & 0 & w & 0 & 0 & -w & 0 \\
   
    0 & 0 & 0 & x & 0 & 0 & 0 & 0 & 0 & 0 & 0 & 0 & 0 & w & 0 & -w \\
    
    z & -z & 0 & z & x & y & -w & w & 0 & -w & 0 & 0 & 0 & 0 & 0 & 0 \\
    
    0 & 0 & 0 & 0 & 0 & 0 & y & -y & 0 & 0 & z & -z & 0 & 0 & 0 & 0  \\
    
    0 & 0 & 0 & 0 & w & 0 & 0 & x & -y & 0 & 0 & 0 & z & -z & 0 & 0  \\
    
    0 & 0 & 0 & 0 & 0 & w & -y & 0 & x & -y & 0 & 0 & 0 & 0 & z & 0   \\
    
    0 & 0 & 0 & 0 & 0 & 0 & 0 & 0 & 0 & x & 0 & 0 & 0 & 0 & 0 & z  \\
    
    0 & 0 & t & 0 & 0 & 0 & 0 & 0 & z & 0 & 0 & x & 0 & y & 0 & 0  \\
    
    0 & 0 & 0 & 0 & 0 & 0 & 0 & 0 & 0 & 0 & -x & 0 & y & 0 & x & -y 
    \end{array}\right)
   \end{equation*}
Since $I$ has a linear presentation and $\rank(\psi)$ attains its maximal possible value, we apply the Jacobian dual criterion and we conclude that $\varphi$ is a birational map onto its image.
\endproof

\begin{remark}\label{rem: B4 and F4 self companion}
It follows from the equation \eqref{eq: unexpected cone for F4} of the unexpected cone that $Y= im(\varphi ) \subset \P^{11}$ is selfdual under the BMSS duality. This turns out to be  true in a much more general set up.  Indeed, let $Z\subset \P^n $  be a finite set of points and let $S_P\subset P^n$ be an unexpected hypersurface of degree $d$ and multiplicity $d$ at a general point $P=(a_0:\cdots :a_n)\in \P^n$ passing through all points in $Z$. Consider $$0=F_Z((a_0:\cdots:a_n),(x_0:\cdots:x_n))=\sum _{i=0}^N
q_{i}(a_0:\cdots:a_n)f_{i}(x_0:\cdots:x_n)$$ the equation of the unexpected cone $S_P$. Set $X=im(\varphi )$  being $ \varphi : \P^n \dasharrow \P^N$ the rational  map defined
 by $f_0,\ldots ,f_N$. By \cite[Theorem 8]{Matrixwise}, $X$ is selfdual under the BMSS duality.
\end{remark}
\subsection{The \texorpdfstring{$H_3$}{H3} root system}
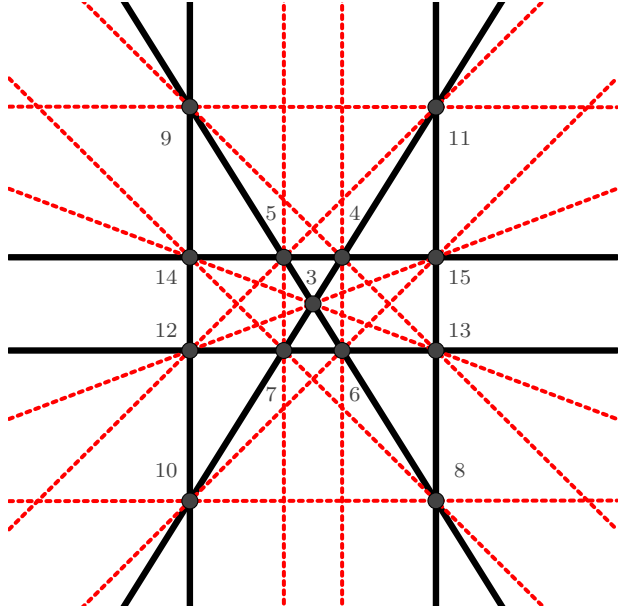
\begin{figure}[h!]
    \centering
\definecolor{qqccqq}{rgb}{1,0,0}
\definecolor{ffqqtt}{rgb}{1,0,0}
\definecolor{ffqqqq}{rgb}{1,0,0}
\definecolor{ttttff}{rgb}{0,0,0}
\begin{tikzpicture}[line cap=round,line join=round,>=triangle 45,x=1cm,y=1cm]
\clip(-4,-4) rectangle (4,4);
\draw [line width=2.4pt,domain=-4:4] plot(\x,{(--0.618-0*\x)/1});
\draw [line width=2.4pt,domain=-4:4] plot(\x,{(-0.618-0*\x)/1});
\draw [line width=2.4pt] (1.618,-4.588760330578508) -- (1.618,4.487122464312542);
\draw [line width=2.4pt] (-1.62,-4.588760330578508) -- (-1.62,4.487122464312542);
\draw [line width=2.4pt,domain=-4:4] plot(\x,{(-0-1*\x)/-0.62});
\draw [line width=2.4pt,domain=-4:4] plot(\x,{(-0-1*\x)/0.62});
\draw [line width=1.5pt,dotted,color=ffqqqq,domain=-4:4] plot(\x,{(--8.455354838709678-0.003225806451613078*\x)/3.238});
\draw [line width=1.5pt,dotted,color=ffqqqq,domain=-4:4] plot(\x,{(-8.455354838709678--0.003225806451613078*\x)/3.238});
\draw [line width=1.5pt,dotted,color=ffqqqq] (-0.38316,-4.588760330578508) -- (-0.38316,4.487122464312542);
\draw [line width=1.5pt,dotted,color=ffqqqq] (0.38316,-4.588760330578508) -- (0.38316,4.487122464312542);
\draw [line width=1.5pt,dotted,color=ffqqqq,domain=-4:4] plot(\x,{(-3.22651741935484--3.230903225806452*\x)/3.238});
\draw [line width=1.5pt,dotted,color=ffqqqq,domain=-4:4] plot(\x,{(--3.2277534193548396--3.2276774193548388*\x)/-3.238});
\draw [line width=1.5pt,dotted,color=ffqqqq,domain=-4:4] plot(\x,{(--3.22651741935484-3.230903225806452*\x)/3.238});
\draw [line width=1.5pt,dotted,color=ffqqqq,domain=-4:4] plot(\x,{(-3.2277534193548396-3.2276774193548388*\x)/-3.238});
\draw [line width=1.5pt,dotted,color=ffqqqq,domain=-4:4] plot(\x,{(-0.0012360000000000149-1.236*\x)/3.238});
\draw [line width=1.5pt,dotted,color=ffqqqq,domain=-4:4] plot(\x,{(--0.0012360000000000149--1.236*\x)/3.238});
\begin{scriptsize}
\draw [fill=uuuuuu] (0,0) circle (3pt);
\draw[color=uuuuuu] (-0.02,0.36) node {$3$};
\draw [fill=uuuuuu] (0.38316,0.618) circle (3pt);
\draw[color=uuuuuu] (0.55,1.2) node {$4$};
\draw [fill=uuuuuu] (-0.38316,0.618) circle (3pt);
\draw[color=uuuuuu] (-0.55,1.2) node {$5$};
\draw [fill=uuuuuu] (0.38316,-0.618) circle (3pt);
\draw[color=uuuuuu] (0.55,-1.2) node {$6$};
\draw [fill=uuuuuu] (-0.38316,-0.618) circle (3pt);
\draw[color=uuuuuu] (-0.55,-1.2) node {$7$};
\draw [fill=uuuuuu] (1.618,-2.609677419354839) circle (3pt);
\draw[color=uuuuuu] (1.93,-2.2) node {$8$};
\draw [fill=uuuuuu] (-1.62,2.612903225806452) circle (3pt);
\draw[color=uuuuuu] (-1.93,2.2) node {$9$};
\draw [fill=uuuuuu] (-1.62,-2.612903225806452) circle (3pt);
\draw[color=uuuuuu] (-1.93,-2.2) node {$10$};
\draw [fill=uuuuuu] (1.618,2.609677419354839) circle (3pt);
\draw[color=uuuuuu] (1.93,2.2) node {$11$};
\draw [fill=uuuuuu] (-1.62,-0.618) circle (3pt);
\draw[color=uuuuuu] (-1.93,-0.35) node {$12$};
\draw [fill=uuuuuu] (1.618,-0.618) circle (3pt);
\draw[color=uuuuuu] (1.93,-0.35) node {$13$};
\draw [fill=uuuuuu] (-1.62,0.618) circle (3pt);
\draw[color=uuuuuu] (-1.93,0.35) node {$14$};
\draw [fill=uuuuuu] (1.618,0.618) circle (3pt);
\draw[color=uuuuuu] (1.93,0.35) node {$15$};
\end{scriptsize}
\end{tikzpicture}
    \caption{The $H_3$ configuration of points with $P_{1}$ and $P_{2}$ at the infinity}
    \label{fig: H3}
\end{figure}

   Associated to the root system $H_3$ we have a set $Z(H_3)$ of  $15$ points, which can be assigned the following coordinates:
$$\begin{array}{lll}
P_{1} = [1:0:0] & P_{2} = [0:1:0] & P_{3} = [0:0:1]
\end{array}$$
$$\begin{array}{llll}
P_{4} = [1:u:u^2] & P_{5} = [-1:u:u^2] & P_{6} = [1:-u:u^2] & P_{7} = [1:u:-u^2] \\
\end{array}$$
$$\begin{array}{llll}
P_{8} = [u:-u^2:1] & P_{9} = [-u:u^2:1] & P_{10} = [u:u^2:-1] & P_{11} = [u:u^2:1] \\
\end{array}$$
$$\begin{array}{llll}
P_{12} = [u^2:1:-u] & P_{13} = [u^2:-1:u] & P_{14} = [-u^2:1:u] & P_{15} = [u^2:1:u] \\
\end{array}$$
where $u^2-u-1=0$, so that $u$ is the golden ratio or its Galois conjugate $\bar{u}=1-u$. The configuration together with $6$ lines, each of them containing $5$ configuration points ($2$ of them are at the infinity) is visualised at Figure \ref{fig: H3}. The $10$ thin lines contain $3$ configuration points each. In order to increase the readability of the figure, only numbers of points are indicated.

The lines with $5$ configuration points have rather simple equations:
$$L_1:\; y-(u-1)z=0,\;\; L_2:\; y+(u-1)z=0,\;\; L_3:\; x-uz=0,\;\; L_4:\; x+uz=0,$$
$$L_5:\; x-(u-1)y=0,\;\; L_6:\; x+(u-1)y=0.$$
The saturated ideal $I(H_3)$ is generated by six quintics. One can choose as the generators the products of all but one line:
\begin{equation}\label{eq: f1..f6}
\begin{aligned}
& f_1=L_2L_3L_4L_5L_6,\;
f_2=L_1L_3L_4L_5L_6,\;
f_3=L_1L_2L_4L_5L_6,\;\\
&
f_4=L_1L_2L_3L_5L_6,\;
f_5=L_1L_2L_3L_4L_6,\;
f_6=L_1L_2L_3L_4L_5.\;
\end{aligned}
\end{equation}
This follows from the fact that the lines $L_1,\ldots,L_6$ form the so called star-configuration, see \cite[Proposition 2.9]{GHM13}, see also \cite[Lemma 3.1 and 3.2]{Schenck11} for an alternative argument. Below we present an equivalent direct argument based on the Hilbert-Burch theorem.
 
The ideal $I(H_3)$ has 13 generators, all in degree 6. Indeed,  $I(H_3)$ is generated by the 6 maximal minors of the $6\times 5$ matrix $U$ with linear entries:
\begin{equation*}
  U =\left(\begin{array}{ccccc}
L_1 & 0 & 0 & 0 & 0 \\
-L_2 & L_2 & 0 & 0 & 0 \\
0 & -L_3 & L_3 & 0 & 0 \\
0 & 0 & -L_4 & L_4 & 0 \\
0 & 0 & 0 & -L_5  & L_5 \\
0 & 0 & 0 & 0 & -L_6
\end{array}\right).
   \end{equation*}
Therefore, $I(H_3)$ has a minimal free resolution:
\begin{equation}\label{exact}
0 \longrightarrow  S(-6)^{5} \xrightarrow{U} S(-5)^{6} \longrightarrow I(H_3) \longrightarrow 0
\end{equation}

\noindent with $S=\C[x,y,z]$ and we conclude that $\dim I(H_3)_6=13$. It is convenient to write down explicitly generators of $I(H_3)_6$:
$$g_i=L_1 f_i, \mbox{ for }i=1,\ldots,6,$$
$$g_7=L_2 f_1,\; g_8=L_2f_3,\; g_9=L_2f_4,\;g_{10}=L_2f_5,\;g_{11}=L_2f_6,\;g_{12}=L_3f_1,\; g_{13}=L_3f_2.$$

The set of points $Z(H_3)$ admits a unique unexpected curve $C\subset \P^2$ of degree 6 and multiplicity 5 at a general point $P=(a:b:c)\in \P^2$ defined by the following equation:
$$
F_{Z(H_3)}((a:b:c),(x:y:z))=\sum_{i=1}^{13}h_i(a:b:c)g_i(x:y:z),
$$
where
$$
h_1=\left(\frac72 u+2\right)\cdot b\cdot\left(b^4+(-44u+72)ab^2c+(-6u+8)b^2c^2+(-4u+12)ac^3+(u-3)c^4\right),
$$
$$
h_2=\frac{-10u-5}{4}\left(a+\frac{u+2}{5}b+\frac{-3u-1}{5}c\right)\left(b+(u-1)c\right)^4,
$$$$
h_3=\frac{2u-1}{4}\left(a+(u+2)b+(u-1)c\right)\left(a-uc\right)^4,
$$
$$
h_4=\frac{-2u+1}{4}\left(a-(u+2)b-(u-1)c\right)\left(a+uc\right)^4,
$$$$
h_5=\frac{4u+3}{4}\left(a-(u-1)b\right)^4\left(a+ub-(u-3)c\right),
$$
$$
h_6=\frac{-4u-3}{4}\left(a+(u-1)b\right)^4\left(a-ub+(u-3)c\right),
$$
$$
h_7=\frac{10u+5}{4}\left(a+\frac{-u-2}{5}b+\frac{-3u-1}{5}c\right)\left(b-(u-1)c\right)^4,
$$$$
h_8=\frac{2u-1}{4}\left(a-(u+2)b+(u-1)c\right)\left(a-uc\right)^4,
$$
$$
h_9=\frac{-2u+1}{4}\left(a+(u+2)b-(u-1)c\right)\left(a+uc\right)^4,
$$
$$
h_{10}=\frac{-4u-3}{4}\left(a-(u-1)b\right)^4\left(a+ub+(u-3)c\right),
$$$$
h_{11}=\frac{4u+3}{4}\left(a+(u-1)b\right)^4\left(a-ub-(u-3)c\right),
$$
$$
h_{12}=\frac{-3u-1}{2}\left(b+uc\right)\left(b-(u-1)c\right)^4,
$$
$$
h_{13}=\frac{3u+1}{2}\left(b-uc\right)\left(b+(u-1)c\right)^4.
$$

\vskip 2mm
\begin{figure}[h!]
    \centering
    \includegraphics[width=0.5\textwidth]{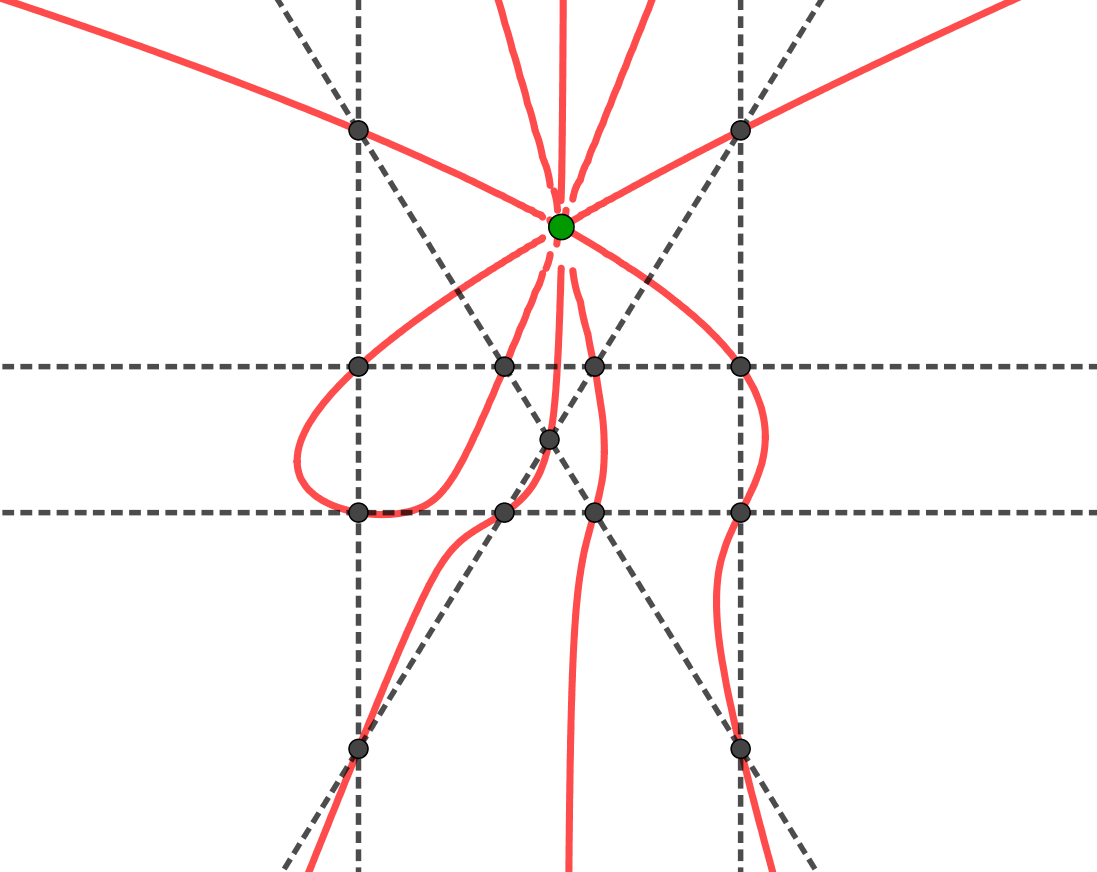}
    \caption{Unexpected sextic with a point of multiplicity $5$}
    \label{fig: H3 curve}
\end{figure}
\noindent    Figure \ref{fig: H3 curve} indicates the unexpected sextic curve with multiplicity $5$ at a general point. The configuration points of the $H_3$ root system are arranged as in Figure \ref{fig: H3}. The lines through only $3$ configuration points and the names of the points are omitted for transparency.

Let $\varphi:\P^2\dashrightarrow\P^{12}$ be the rational map defined by the 13 generators $g_1,\ldots ,g_{13}$ of $I(Z(H_3))$.
\begin{proposition}\label{prop: surface of segree 21}
   The image of $\varphi$ is an arithmetically Cohen-Macaulay smooth surface $X\subset \P^{12}$ of degree $21$ whose homogeneous ideal $I(X)$ is generated by quadrics.
\end{proposition}   
\proof
   We consider $\pi :Y \to\P^2$  the blow up of $\P^2$ at the  $15$ points of $Z(H_3)$. Set $H=\pi ^*\calo_{\P^2}(1)$ and denote by $E_i$, $i=1,\ldots ,15,$ the exceptional divisors.
   
   \noindent {\bf Claim 1:} The linear system $L=6H-\sum _{i=1}^{15}E_i$ is base point free and very ample.
   
   \noindent {\em Proof of Claim 1.} We sheafify the exact sequence \eqref{exact} and we get the exact sequence:
   $$
   0 \longrightarrow  \calo_{\P^2}(-6)^5 \longrightarrow  \calo_{\P^2}(-5)^6 \longrightarrow  {\cal I} _{Z(H_3)} \longrightarrow 0 .
   $$
We compute the cohomology and we get $H^1 (\mathcal{I}_{Z(H_{3})}(5))=0$. Therefore, the linear system $L=6H-\sum _{i=1}^{15}E_i$ is base point free (see \cite[Theorem 1.4]{DG}). Moreover, since no line in $\P^2$ contains 6 points of $Z(H_3)$, we can apply \cite[Theorem 3.1]{DG} and we conclude that $L$ is very ample.

Since  $L=6H-\sum _{i=1}^{15}E_i$ is base point free and very ample, it defines an embedding, i.e. an isomorphism  onto its image 
    $\varphi_{L}:Y\to\P^{12}$ which lifts the rational map
    $\varphi:\P^2\dashrightarrow\P^{12}.$
 Therefore, we have the commutative diagram
 $$ \begin{tikzcd} 
 Y \arrow[d,"\pi "] \arrow[dr,"\varphi _{L}"] &\\
 \P^2 \arrow[r ,dashrightarrow, "\varphi"] & \P^{12}.
 \end{tikzcd}
 $$
 and the image of $\varphi $ is a smooth surface $X\subset \P^{12}$ of degree $L^2=21$.
 
  Let us now prove that $X$ is an arithmetically Cohen-Macaulay surface.  First we compute the $h$-vector $(h_{Z(H_3)}(i))$ of the set of points $Z(H_3)\subset \P^2$ where 
   $$h_{Z(H_3)}(i):=\Delta H(Z(H_3),i)=H(Z(H_3),i)-H(Z(H_3),i-1)$$
    $$=\dim\left(\frac{k[x,y,z]}{I(Z(H_3))}\right)_i - \dim\left(\frac{k[x,y,z]}{I(Z(H_3))}\right)_{i-1}.$$
   The $h$-vector of $Z(H_3)$ is $(1, 2, 3, 4, 5)$. Therefore, $h_{Z(H_3)}(i)=0$ for all $i\ge 5$ and $\sigma =5$, where $\sigma $ is the least integer $t$ such that $\Delta H(Z(H_3),t)=0$. By \cite[Theorem B]{GG}, $X$ is an arithmetically Cohen-Macaulay surface and by \cite[Theorem C]{GG}, $I(X)$ is generated by quadrics which proves what we want.
\endproof
  \noindent Turning to the companion surface we take a closer look at the polynomials $h_1,h_2,\ldots,h_{13}.$ First of all, it can be checked directly that $h_1$ is a linear combination of the remaining $12$
   polynomials. These polynomials, in turn, can be paired according to their fourth power factors. Removing the constant terms and taking linear combinations, we obtain the following $12$ quintics in variables $(a:b:c)$:
$$q_1=a\,(b-(u-1)c)^4,\;\; 
  q_2=(b+uc)(b-(u-1)c)^4,$$   
$$q_3=a\,(b+(u-1)c)^4,\;\;
  q_4=(b-uc)(b+(u-1)c)^4,$$  
$$q_5=b\,(a-uc)^4,\;\;
  q_6=(a+(u-1)c)(a-uc)^4,$$  
$$q_7=b\,(a+uc)^4,\;\;
  q_8=(a-(u-1)c)(a+uc)^4,$$  
$$q_9=c\,(a-(u-1)b)^4,\;\;
  q_{10}=(a+ub)(a-(u-1)b)^4,$$  
$$q_{11}=c\,(a+(u-1)b)^4,\;\;
  q_{12}=(a-ub)(a+(u-1)b)^4.$$  
Interestingly, the unexpected curve written down in these generators has the following form
$$\begin{aligned}
F_{Z(H_3)}& =\frac{1}{2}\Bigl[ 
 (3u+1)\;\left[M_1f_1\; q_1-xf_1\; q_2-M_2f_2\; q_3+xf_2\; q_4\right]\\
&+ (1-2u)\;\left[M_3f_3\; q_5-\;yf_3\; q_6-\;M_4f_4\; q_7+\;yf_4\; q_8\right]\\
&+ (3u+1)\;\left[M_5f_5\; q_9-\;zf_5\; q_{10}-\;M_6f_6\; q_{11}+\;zf_6\; q_{12}\right]\Bigr],
\end{aligned}$$
where $M_1,\ldots,M_6$ are Galois conjugates of lines $L_1,\ldots,L_6$ and $f_1,\ldots,f_6$ are defined in \eqref{eq: f1..f6}:
$$
\begin{aligned}
&M_1: \; y+uz=0,\;\; M_2: \; y-uz=0, \\
&M_3: \; x-\bar uz=0, \;\; M_4: \; x+\bar uz=0,\\
&M_5:\; x+uy=0,\;\; M_6:\; x-uy=0 .
\end{aligned}
$$
These lines are indicated as dashed lines in Figure \ref{fig: conj H3}, the $L$-lines are indicated there as the thin lines. The dashed lines intersect in the Galois conjugate $H_3$ configuration of points. In particular the origin and the two points at the infinity belong to both configurations.
\begin{figure}[h!]
    \centering
\definecolor{qqccqq}{rgb}{1,0,0}
\definecolor{ffqqtt}{rgb}{0,1,0}
\definecolor{ffqqqq}{rgb}{1,0,0}
\definecolor{ttttff}{rgb}{0,0,0}
\begin{tikzpicture}[line cap=round,line join=round,>=triangle 45,x=1cm,y=1cm]
\clip(-4,-4) rectangle (4,4);
\draw [line width=1pt,color=ffqqtt,domain=-4:4] plot(\x,{(--0.618-0*\x)/1});
\draw [line width=1pt,color=ffqqtt,domain=-4:4] plot(\x,{(-0.618-0*\x)/1});
\draw [line width=1pt,color=ffqqtt] (1.618,-4.588760330578508) -- (1.618,4.487122464312542);
\draw [line width=1pt,color=ffqqtt] (-1.62,-4.588760330578508) -- (-1.62,4.487122464312542);
\draw [line width=1pt,color=ffqqtt,domain=-4:4] plot(\x,{(-0-1*\x)/-0.62});
\draw [line width=1pt,color=ffqqtt,domain=-4:4] plot(\x,{(-0-1*\x)/0.62});
\draw [line width=2.4pt,dashed,domain=-4:4] plot(\x,1.618);
\draw (3.5,1.3) node {$M_1$};
\draw [line width=2.4pt,dashed,domain=-4:4] plot(\x,-1.618);
\draw (3.5,-1.3) node {$M_2$};
\draw [line width=2.4pt,dashed] (0.618,-4) -- (0.618,4);
\draw (1,3) node {$M_4$};
\draw [line width=2.4pt,dashed] (-0.618,-4) -- (-0.618,4);
\draw (-1,3) node {$M_3$};
\draw [line width=2.4pt,dashed,domain=-4:4] plot(\x,{-0.618*\x});
\draw (-3.5,2.7) node {$M_5$};
\draw [line width=2.4pt,dashed,domain=-4:4] plot(\x,{0.618*\x});
\draw (-3.5,-2.7) node {$M_6$};
\begin{scriptsize}
\draw [fill=uuuuuu] (0,0) circle (2pt);
\draw [fill=uuuuuu] (0.38316,0.618) circle (2pt);
\draw [fill=uuuuuu] (-0.38316,0.618) circle (2pt);
\draw [fill=uuuuuu] (0.38316,-0.618) circle (2pt);
\draw [fill=uuuuuu] (-0.38316,-0.618) circle (2pt);
\draw [fill=uuuuuu] (1.618,-2.609677419354839) circle (2pt);
\draw [fill=uuuuuu] (-1.62,2.612903225806452) circle (2pt);
\draw [fill=uuuuuu] (-1.62,-2.612903225806452) circle (2pt);
\draw [fill=uuuuuu] (1.618,2.609677419354839) circle (2pt);
\draw [fill=uuuuuu] (-1.62,-0.618) circle (2pt);
\draw [fill=uuuuuu] (1.618,-0.618) circle (2pt);
\draw [fill=uuuuuu] (-1.62,0.618) circle (2pt);
\draw [fill=uuuuuu] (1.618,0.618) circle (2pt);
\end{scriptsize}
\end{tikzpicture}
    \caption{The $H_3$ configuration of points and the Galois conjugate lines}
    \label{fig: conj H3}
\end{figure}
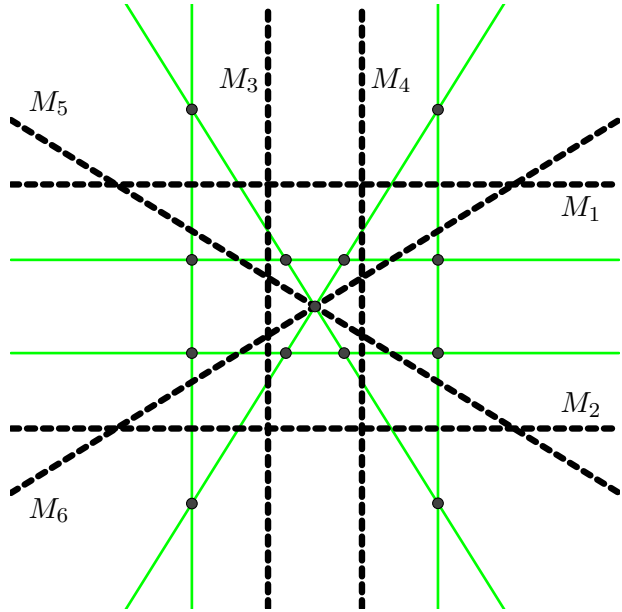 
\begin{remark}\label{rem: H3 companion}
  We describe the companion surface $X'$ of $X$.
  
  Let $ \psi: \P^2 \dashrightarrow \P^{11}$ be the rational map defined by the 12 generators $q_1,\ldots,q_{12}$ and we call $X'$ the image of $\psi $. To check that the map $\psi:\P^2 \dashrightarrow  \P^{11}$ has degree one and even more it is a birational map onto its image we use the packages \texttt{"Cremona"} and
\texttt{"RationalMaps"} of \textit{Macaulay2} as before.
  
    \vskip 4mm
\noindent\texttt{loadPackage "Cremona"}\\
\texttt{loadPackage "RationalMaps"}\\
\texttt{kk = ZZ/22621}\\
\texttt{R = kk[a,b,c]}\\
\texttt{u = 1873 - - golden ratio}\\
\texttt{I =  ideal(a\,(b-(u-1)c)^4, (b+uc)(b-(u-1)c)^4, a\,(b+(u-1)c)^4, (b-uc)(b+(u-1)c)^4,}\\ 
\texttt{b\,(a-uc)^4,(a+(u-1)c)(a-uc)^4, b\,(a+uc)^4,(a-(u-1)c)(a+uc)^4, c\,(a-(u-1)b)^4,}\\
\texttt {(a+ub)(a-(u-1)b)^4,  c\,(a+(u-1)b)^4,(a-ub)(a+(u-1)b)^4)}\\
\texttt{degree rationalMap gens I}\\    
(Using the package "Cremona" we get that the degree is 1)\\
\texttt{S=kk[q\_1..q\_(numgens I)]}\\
\texttt{varphi = map(R, S, gens I)}\\
\texttt{isBirationalOntoImage(varphi) - - true}\\
(Using the package "RationalMaps" we get that the map is birational onto its image)
   \vskip 4mm

      Once we know that $\psi$ is a birational map onto its image we compute the degree of $X'$, it is $25$. The smoothness of $X'$ follows from the Jacobian criterion and a  computation which shows that  the rank of $\begin{pmatrix} \frac{\partial q_{i}}{\partial a}& \frac{\partial q_{i}} {\partial b} & \frac{\partial q_{i}} {\partial c }\end{pmatrix} _{i=1,\ldots ,12}$ is 3. It is a projection of the Veronese $V_{2,5}$ by a subspace of dimension eight that does not intersect the Secant variety. Moreover, the ideal of $X'$ is generated by  quadrics and  cubics.
  
\end{remark}
 \begin{remark}\label{rem: H3 bicompanion}
 We consider the rational map $\bar \varphi: \mathbb P ^2 \dashrightarrow \mathbb P^{11} $ where
 $$\bar \varphi=(M_1f_1,xf_1,M_2f_2,xf_2,M_3f_3,yf_3,M_4f_4,yf_4,M_5f_5,zf_5,M_6f_6,zf_6) $$
  and we call $\bar X$ its image. Using  {\em Macaulay2} we get that $\bar \varphi $ is a birational map onto  $\bar X$, $degree(\bar X)=21$ and $I(\bar X)$ is generated by 32 quadrics. The surface $\bar X\subset \P ^{11}$ can be seen as a projection of the surface $X\subset \P ^{12}$ described in Proposition 4.5. However, $\bar X$ is not aCM, as its $h$-vector is $h:(1, \; 9,\; 13,\; -3,\; 1)$, so it is not positive.
 \end{remark}
\section{Companion surfaces for Fermat arrangements of lines}
\label{sec: Fermats}
   Recall that in suitable coordinates the points in the root system $B_3$ are dual to linear factors of the Fermat-type polynomial
   $$F_{2,3}(x:y:z)=xyz(x^2-y^2)(y^2-z^2)(z^2-x^2).$$
   Thus it is natural to generalize $B_3$ from the perspective of Fermat-type arrangements. We follow this path of thoughts in this section. Our notation is consistent with the notation introduced in \cite{Szp19c}.
   
   Let
   $$F_{m,0}(x:y:z)=(x^m-y^m)(y^m-z^m)(z^m-x^m)$$
   and let $Z(F_{m,0})$ be the set of points dual to the linear factors of $F_{m,0}$.
   Thus $Z(F_{m,0})$ is a set of $3m$ points distributed evenly on the coordinate lines.
   This is visualised in Figure \ref{fig: Z(F(5,0))} for $m=5$.
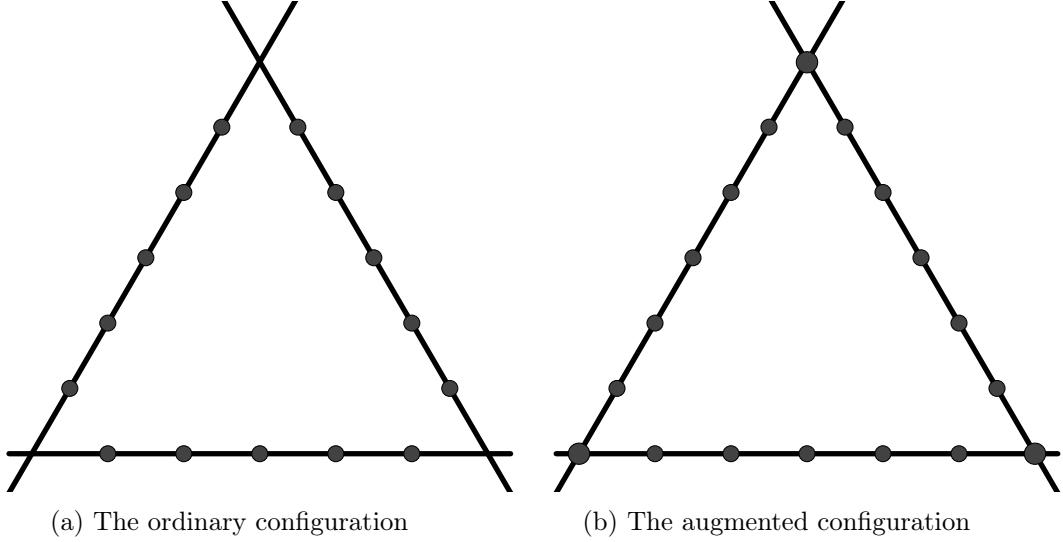
\begin{figure}[h!]
\centering
\begin{subfigure}{.48\textwidth}
  \centering
  \begin{tikzpicture}[line cap=round,line join=round,>=triangle 45,x=1cm,y=1cm]
\clip(-1,-.5) rectangle (7,6);
\draw [line width=2pt,domain=-4.94:14.18] plot(\x,{(-0--5.196152422706633*\x)/3});
\draw [line width=2pt,domain=-4.94:14.18] plot(\x,{(--31.1769145362398-5.196152422706633*\x)/3});
\draw [line width=2pt,domain=-.3:6.3] plot(\x,{(-0-0*\x)/-6});
\begin{scriptsize}
\draw [fill=uuuuuu] (0.5,0.8660254037844393) circle (3pt);
\draw [fill=uuuuuu] (1,1.7320508075688785) circle (3pt);
\draw [fill=uuuuuu] (1.5,2.5980762113533182) circle (3pt);
\draw [fill=uuuuuu] (2,3.4641016151377584) circle (3pt);
\draw [fill=uuuuuu] (2.5,4.330127018922198) circle (3pt);
\draw [fill=uuuuuu] (3.5,4.330127018922196) circle (3pt);
\draw [fill=uuuuuu] (4,3.464101615137759) circle (3pt);
\draw [fill=uuuuuu] (4.5,2.5980762113533244) circle (3pt);
\draw [fill=uuuuuu] (5,1.7320508075688847) circle (3pt);
\draw [fill=uuuuuu] (5.5,0.8660254037844328) circle (3pt);
\draw [fill=uuuuuu] (1,0) circle (3pt);
\draw [fill=uuuuuu] (2,0) circle (3pt);
\draw [fill=uuuuuu] (3,0) circle (3pt);
\draw [fill=uuuuuu] (4,0) circle (3pt);
\draw [fill=uuuuuu] (5,0) circle (3pt);
\end{scriptsize}
\end{tikzpicture}
  \caption{The ordinary configuration}
  \label{fig: Z(F(5,0))}
\end{subfigure}%
\begin{subfigure}{.48\textwidth}
  \centering
    \begin{tikzpicture}[line cap=round,line join=round,>=triangle 45,x=1cm,y=1cm]
\clip(-1,-.5) rectangle (7,6);
\draw [line width=2pt,domain=-4.94:14.18] plot(\x,{(-0--5.196152422706633*\x)/3});
\draw [line width=2pt,domain=-4.94:14.18] plot(\x,{(--31.1769145362398-5.196152422706633*\x)/3});
\draw [line width=2pt,domain=-.3:6.3] plot(\x,{(-0-0*\x)/-6});
\begin{scriptsize}
\draw [fill=uuuuuu] (0.5,0.8660254037844393) circle (3pt);
\draw [fill=uuuuuu] (1,1.7320508075688785) circle (3pt);
\draw [fill=uuuuuu] (1.5,2.5980762113533182) circle (3pt);
\draw [fill=uuuuuu] (2,3.4641016151377584) circle (3pt);
\draw [fill=uuuuuu] (2.5,4.330127018922198) circle (3pt);
\draw [fill=uuuuuu] (3.5,4.330127018922196) circle (3pt);
\draw [fill=uuuuuu] (4,3.464101615137759) circle (3pt);
\draw [fill=uuuuuu] (4.5,2.5980762113533244) circle (3pt);
\draw [fill=uuuuuu] (5,1.7320508075688847) circle (3pt);
\draw [fill=uuuuuu] (5.5,0.8660254037844328) circle (3pt);
\draw [fill=uuuuuu] (0,0) circle (4pt);
\draw [fill=uuuuuu] (1,0) circle (3pt);
\draw [fill=uuuuuu] (2,0) circle (3pt);
\draw [fill=uuuuuu] (3,0) circle (3pt);
\draw [fill=uuuuuu] (4,0) circle (3pt);
\draw [fill=uuuuuu] (5,0) circle (3pt);
\draw [fill=uuuuuu] (6,0) circle (4pt);
\draw [fill=uuuuuu] (3,5.19) circle (4pt);
\end{scriptsize}
\end{tikzpicture}
  \caption{The augmented configuration}
  \label{fig: Z(F(5,3))}
\end{subfigure}
\caption{Fermat configurations of points for $m=5$}
\label{fig:test}
\end{figure}

   We recall a result \cite[Theorem 4.6]{Szp19c} which motivated our present research.
\begin{theorem}[Szpond]\label{thm: Fermat unexpected}
   Let $m\geq 5$ be an integer and let $R=(a:b:c)$ be a general point in $\P^2$.
   The set $Z(F_{m.0})$ admits a unique irreducible unexpected curve $C_{R,m}$
   of degree $m+2$ and multiplicity $m+1$. Moreover, all these curves $C_{R,m}$
   pass through the coordinate points.
\end{theorem}
Set  $F_{m,3}=xyzF_{m,0}$.  It turns out that for the augmented configuration
   $$Z(F_{m,3})=Z(F_{m,0})\cup(1:0:0)\cup(0:1:0)\cup(0:0:1)$$
   the curve $C_{R,m}$ is unexpected as soon as $m\geq 2$.
   For $m=2$ we get exactly the unexpected curve associated to the $B_3$ root system.
   
\begin{proposition}\label{prop: ideal of Z(F(m,3)) m even}
   The ideal $I(Z(F_{m,3}))\subset S=\C[x,y,z]$ is generated by
   $$xyz,\; G_x=yz(y^m-(-z)^m),\; G_y=zx(z^m-(-x)^m),\; G_z=xy(x^m-(-y)^m).$$
   Even more,  $I(Z(F_{m,3}))$ is  defined by the maximal minors of the matrix
   $$M=\begin{pmatrix} y^m-(-z)^m & z^m-(-x)^m & x^m-(-y)^m \\
    x & 0 & 0 \\
    0 & y & 0 \\
    0 & 0 & z 
    \end{pmatrix}
   $$
   and, hence, it has  a minimal free resolution of the following type:
   $$
   0 \longrightarrow S(-m-3)^3 \xrightarrow{M} S(-3)\oplus S(-m-2)^3 \longrightarrow  I(Z(F_{m,3})) \longrightarrow 0.
   $$
\end{proposition}   
\proof
   First we note that the reason for the sign in the three binomial generators depending on the parity of $m$ is due to the convention that a point $(p:q:r)$ is dual to the line $px+qy+rz=0$. Thus, for example, if $m=3$ and $\varepsilon$ is a primitive root of the unity of order $3$, then 
   $$(x-\varepsilon y)(x-\varepsilon^2 y)(x-y)=x^3-y^3,$$
   and the dual points are
   $$(1:-\varepsilon:0),\; (1:-\varepsilon^2:0),\; (1:-1:0).$$
   Their ideal is generated by $(z,\; x^3+y^3)$.
   
   It is now easy to see that for arbitrary $m$ the trace of $Z(F_{m,0})$ on the coordinate line, 
   say $x=0$, is a set of zeroes of a complete intersection ideal generated by 
   $$x\;\mbox{ and }\; yz(y^m+(-1)^{m+1} z^m)$$
   and similarly on the other lines. Thus
   \begin{align*}
      I(Z(F_{m,3})) = & (xy,yz,zx)\cap(x,yz(y^m+(-1)^{m+1} z^m))\\
      & \cap(y,xz(z^m+(-1)^{m+1} x^m))\cap(z,xy(x^m+(-1)^{m+1} y^m)).
   \end{align*}   
   and the claim follows readily. 
\endproof
We are interested in homogeneous elements $[I(Z(F_{m,3}))]_{m+2}$
 in $I(Z(F_{m,3}))$ of degree $(m+2)$.
   They are $G_x, G_y, G_z$ and monomials of the form
  $$xyz \lambda\mbox{ where } \lambda \mbox{ is a monomial of degree } m-1.$$ Hence
\begin{equation}\label{eq: dim in deg m+2}
   \dim [I(Z(F_{m,3}))]_{m+2}=\frac{1}{2}m^2+\frac{1}{2}m+3.
\end{equation}
   Geometrically the elements in $[I(Z(F_{m,3}))]_{m+2}$ correspond to elements of the linear system of forms of degree $m+2$ vanishing at all points of the set $Z(F_{m,3})$.
\subsection{Positivity on anticanonical surfaces}
\label{sec: positivity}
   Before we continue the study of the linear system determined by $[I(Z(F_{m,3}))]_{m+2}$, we need to recall some results
   on anticanonical systems on surfaces. Results presented here were obtained by Harbourne in \cite{Har97JA} and \cite{Har97TAMS}. 
   We begin by the following useful Lemma, which is a combination of Lemma II.2 and Corollary II.3 in \cite{Har97TAMS} and Lemma 2.2 and Corollary 2.3 in \cite{Har97JA}.
\begin{lemma}[Harbourne]\label{lem: nef on rational}
   Let $Y$ be a smooth rational surface and let $D$ be a nef class on $Y$. Then
$$
   h^2(Y,D)=0\;\;\mbox{ and }\;\; D^2\geq 0.
$$
   Moreover, if $-K_Y$ is effective, then so is $D$.
\end{lemma}   
   Assuming the effectivity of the anticanonical class, one can say in fact much more. The following statement is extracted from a much more precise (and complicated) result of Harbourne \cite[Theorem III.1]{Har97TAMS} and \cite[Theorem 2.11]{Har97JA}.
\begin{proposition}[Harbourne]\label{prop: nef on anticanonical}
   Let $Y$ be a smooth anticanonical surface and let $D$ be a nef class on $Y$ such that $-K_Y.D\geq 2$. Then $D$ is non-special, i.e., $h^1(Y,D)=0$ and the linear system $|D|$ is base point free.
\end{proposition}   
   Thus the anticanonical degree at least $2$ implies that $|D|$ defines a morphism. If this degree is at least $3$, then the positivity of $|D|$ increases, see \cite[Proposition 3.2]{Har97JA}.
\begin{proposition}[Harbourne]\label{prop: nef on anticanonical geq 3}
   Let $Y$ be a smooth anticanonical surface and let $D$ be a nef class on $Y$ such that $-K_Y.D\geq 3$ and $D^2>0$. Then the morphism $\varphi_D$ defined by the elements of $|D|$ is birational and its image $Y'$ is a normal surface obtained by contracting all curves perpendicular to $|D|$. Moreover a general member of $|D|$ is smooth and irreducible.
\end{proposition}
\subsection{On the positivity of maps determined by Fermat arrangements} 
   Let $m$ be fixed and let $\eps$ be a primitive root of the unity of order $m$. It is convenient to introduce the following notation:
   $$P_{x,\alpha}=(0:1:\eps^\alpha),\;
   P_{y,\alpha}=(1:0:\eps^\alpha),\;
   P_{z,\alpha}=(1:\eps^\alpha:0)$$
   for $\alpha=0,\ldots,m-1$ and
   $$P_{x,y}=(0:0:1),\;
   P_{x,z}=(0:1:0),\;
   P_{y,z}=(1:0:0).$$
   Our next result concerns the positivity of the linear system corresponding to $[I(Z(F_{m,3}))]_{m+2}.$
\begin{theorem}\label{thm: Fermat Z(Fm,3)}
   Let $f:Y\to\P^2$ be the blow up at the points of the set $Z(F_{m,3})$ with the exceptional divisors $E_{a,\alpha}$, $E_{a,b}$ over the points $P_{a,\alpha}$ and $P_{a,b}$ respectively, where $a,b\in\left\{x,y,z\right\}$ and $\alpha\in\left\{0,\ldots,m-1\right\}$. Let $H=f^*\calo_{\P^2}(1)$ be the pull-back of the hyperplane bundle and let $\E$ be the union of all exceptional divisors of $f$. Then the linear system
   $$L=(m+2)H-\E$$
   is base point free and defines a morphism $\varphi_L$ birational onto its image
   $$\varphi_L:Y\to X\subset\P^N,\;\; (N=\frac{1}{2}(m^2+m)+2),$$
 which is an isomorphism away from the proper transforms of the coordinate lines.
\end{theorem}   
\proof
   Let
\begin{align*}
   N_x &= H-E_{x,y}-E_{x,z}-E_{x,0}-\cdots-E_{x,m-1}\\    
   N_y &= H-E_{x,y}-E_{y,z}-E_{y,0}-\cdots-E_{y,m-1}\\    
   N_z &= H-E_{x,z}-E_{y,z}-E_{z,0}-\cdots-E_{z,m-1}
\end{align*}
   be the proper transformations of the coordinate lines.
   
   Our first claim is that $L$ is nef and the only curves having the intersection number with $L$ equal zero are $N_x, N_y$ and $N_z$. 
   
   It is obvious that $L.N_a=0$ for all $a\in\left\{x,y,z\right\}$, so it is enough to check that $L$ has positive intersection with all other irreducible curves on $Y$. Certainly it is $L.E=1$ for all exceptional divisors of $f$. Let $C$ be an irreducible plane curve of degree $d$, different from the coordinate lines, passing through the points $P_{a,\alpha}$ and $P_{a,b}$ with multiplicities $m_{a,\alpha}$ and $m_{a,b}$ respectively. Let
   $$\wtilde{C}=dH-\sum m_{a,\alpha}E_{a,\alpha}-\sum m_{a,b}E_{a,b}$$
   be the proper transform of $C$. Then
\begin{align}
    0\leq \wtilde{C}.N_x &=d-m_{x,y}-m_{x,z}-m_{x,0}-\cdots-m_{x,m-1}, \label{eq: Nx}\\
    0\leq \wtilde{C}.N_y &=d-m_{x,y}-m_{y,z}-m_{y,0}-\cdots-m_{y,m-1}, \label{eq: Ny}\\
    0\leq \wtilde{C}.N_z &=d-m_{x,z}-m_{y,z}-m_{z,0}-\cdots-m_{z,m-1}.\label{eq: Nz}
\end{align}   
   Adding the inequalities \eqref{eq: Nx}, \eqref{eq: Ny} and \eqref{eq: Nz} and rearranging terms we get
$$
   3d-2(m_{x,y}+m_{x,z}+m_{y,z})-\sum m_{a,\alpha}\geq 0.
$$
   Since
$$
   M.\wtilde{C}=4d-m_{x,y}-m_{x,z}-m_{y,z}-\sum m_{a,\alpha}
$$
   we conclude that $M.\wtilde{C}>0$ as asserted.
   
   For the rest of the proof, the key observation is that the surface $Y$ is an anticanonical surface, i.e.,  $-K_Y$ is an effective divisor. Indeed
$$
   -K_Y=N_x+N_y+N_z+E_{x,y}+E_{x,z}+E_{y,z}.
$$
   Since
$$
   -K_Y.M=3
$$
   Proposition \ref{prop: nef on anticanonical} implies that $|M|$ is non-special and base point free. Revoking Lemma \ref{lem: nef on rational} we have thus the vanishing
$$
   h^1(Y, L)=h^2(Y,D)=0,
$$
   which together with the Riemann-Roch formula on $Y$ reconfirms the computation in \eqref{eq: dim in deg m+2}. 
   Since
$$
   L^2=m^2+m+1>0,
$$
   it follows from Proposition \ref{prop: nef on anticanonical geq 3} that the morphism defined by $|L|$ is birational and its image is a normal surface of degree $m^2+m+1$ in $\P^N$ with $N=\frac{1}{2}(m^2+m)+2$ singular in exactly $3$ points, which are images of the proper transforms of the coordinate lines in $\P^2$.
\endproof
   We now turn to the equations defining $X$. To this end let us be a little bit more specific about the coordinates of the rational map $\varphi$ determined by the commutative diagram
$$ \begin{tikzcd} 
 Y \arrow[d,"f"] \arrow[dr,"\varphi_{M}"] &\\
 \P^2 \arrow[r ,dashrightarrow, "\varphi"] & X\subset\P^{N},
 \end{tikzcd}
$$
  which is not just incidentally similar to that in Proposition \ref{prop: threefold of deg 40}. We can assume that
  $$\varphi=(xyz\mu_d:G_x:G_y:G_z),$$
  where $d=m-1$ and 
  $$\mu_d:\P^2\hookrightarrow\P^{N-3}$$
  denotes the Veronese embedding of degree $d$ and the rational map $xyz\mu_d$ is given by coordinates of $\mu_d$ multiplied by the monomial $xyz$.
  
  It is well known, see for example \cite[Corollary 7.2.3]{Weyman-book}, that ideal of the image $V_d$ of $\mu_d$ is generated by $\frac{(d-1)d(d+1)(d+6)}{8}$ quadratic binomials.
  
  Consider the projection $\pi:\P^N\dashrightarrow \P^{N-3}$ determined by the plane $\Pi$ spanned by the last three coordinate points in $\P^N$. Then $\pi(X)=V_d$ and the mapping $\varphi$ (or equivalently $\varphi_L$) can be considered as the \lq\lq unprojection'' of the Veronese surface in $\P^{N-3}$. We use the word unprojection to indicate a morphism inverse to a projection, not quite in the sense as introduced to the birational geometry by Miles Reid.
\begin{proposition}\label{prop_image_Fm}
   Let $m\geq 4$ be a fixed integer and let $I=I(X)$ be the ideal of $X\subset\P^N$. It holds:
   
   (1) $X\subset \P^N$ is an aCM surface.
   
   (2) The ideal $I$ is generated by forms of degree $\le 3$. More precisely,  $I$ is generated by  $3 \binom{m+2}{4}$ quadrics and one cubic $C_m$, where for $m$ odd
   $$C_m:\;\; 2H_xH_yH_z+G_xH_yH_z+H_xG_yH_z+H_xH_yG_z-G_xG_yG_z,$$
   and for $m$ even
      $$C_m:\;\; G_xH_yH_z+H_xG_yH_z+H_xH_yG_z-G_xG_yG_z,$$
   where $H_x=x^myz, H_y=xy^mz$ and $H_z=xyz^m$.
\end{proposition}  
\proof
(1) Let $L\cong \P^{N-1}\subset \P^N$ be a general hyperplane. Since $X\subset \P^N$ is an aCM surface if and only if a general hyperplane section $C= X\cap L$ of $X$ is an aCM curve, it is enough to check that indeed $C\subset L=\P^{N-1}$ is an aCM curve. We observe that 
$$ \deg(C)=\deg(X)=m^2+m+1.$$
By the adjunction formula we have:
$$ 2g(C)-2=C(C+K_X)=C^2-3$$ 
(because $CK_X=LK_X=((m+2)H-\sum_{i=1}^{3m+3}E_i)(-3H+\sum_{i=1}^{3m+3}E_i)=-3$). Therefore, $C^2=2g(C)+1$. Thus, ${\cal O}_C(C)$ is a very ample line bundle on $C$ and by the Corollary to \cite[Theorem 6]{Mum}, ${\cal O}_C(C)$ is projectively normal. So, ${\cal O}_C(C)$ embeds $C$ in a projective space as a projectively normal curve. Since $X$ is rational, it is regular (i.e. $H^1(X, {\cal O}_X)=0$) and from the exact cohomology sequence associated to the exact sequence
$$
0 \longrightarrow {\cal O}_X \longrightarrow {\cal O}_X (C) \longrightarrow {\cal O}_C (C) \longrightarrow 0
$$
we deduce that the projective space into which ${\cal O}_C(C)$ embeds $C$ is just $L$ and we are done.

\vskip 2mm
(2) Let us first compute the dimension of $I(X)_2.$ By \cite[Corollary 7.2.3]{Weyman-book} there are 
    $\frac{(m-2)(m-1)m(m+5)}{8}$ quadrics involving the monomials $xyz(x,y,z)^{m-1}$. A straightforward computation shows that there are $\binom{m}{2}$ quadrics involving $G_x$ (resp., $G_y$ and $G_z$) and the monomials in $xyz(x,y,z)^{m-1}$. Therefore $\dim I(X)_2=3 \binom{m+2}{4}$.  To prove that $I(X)$ is generated by quadrics and cubics, we consider the $h$-vector of $X$. Recall that the $h$-vector of $X$ is a sequence $(h_X(i))_{i\ge 0}$ of non-negative integers satisfying:
   
    \vskip 2mm
    \begin{itemize}
    \item[(i)] For all $i\ge 0$, $h_X(i)=\Delta ^3H(X,i)$ where $H(X,i)=\dim (\C[y_0,\ldots y_N]/I(X))_{i}$.
    \item[(ii)] $\deg(X)=\sum _{i\ge 0}h_X(i)$.
    \item[(iii)] If $\sigma $ is the least integer such that $h_X(i)=0$, then $I(X)$ is generated by forms of degree $\le \sigma$.
    \end{itemize}
   
   \vskip 2mm
    Since $X$ is non-degenerate and $\dim I(X)_2=3 \binom{m+2}{4}$, we have $h_X(0)=1$, $h_X(1)=h_X(2)=N-2$. From the equality $$\deg(X)=m^2+m+1=\sum _{i\ge 0}h_X(i)=2N-3+\sum _{i\ge 3}h_X(i),$$ we conclude that $h_X(i)=0$ for $i\ge 3$, $\sigma =3$ and $I(X)$ is generated by forms of degree $\le 3$. 
    Let us now prove that $I(X)$ is not generated by quadrics i.e. at least one cubic is required. Since $X$ is an aCM surface, we know that the graded Betti numbers of $X$ and  of a general hyperplane section $C=X\cap \P^{N-1}$ are equal; in particular, the degrees of a minimal set of generators of $C$ and $X$ coincide. So, it will be enough to check that a minimal set of generators of $I(C)$ contains a cubic. In \cite[Theorem 2]{GL}, Green and Lazarsfeld proved that $I(C)$ fails to be generated by quadrics if and only if $C$ is hyperelliptic or ${\cal O}_C (C)$ embeds $C$ with a trisecant line, i.e. $H^0(C; C-K_C)\ne 0$ where $K_C$ denotes the canonical divisor of $C$. We consider the exact sequence  
    $$
0 \longrightarrow {\cal O}_X \longrightarrow {\cal O}_X (C) \longrightarrow {\cal O}_C (C) \longrightarrow 0
$$ and its associated exact cohomology sequence:
$$
0 \longrightarrow H^0(X;-C-K_X) \longrightarrow H^0(X;-K_X) \longrightarrow H^0(C;{\cal O}_C(-K_X)) \longrightarrow H^1(X;-C-K_X)\cdots .
$$

From the adjunction formula  $K_C\cong {\cal O}_C(C+K_X)$, we deduce that  $H^0(C;{\cal O}_C(-K_X))\cong H^0(C;{\cal O}_C(C-K_C))$. On the other hand, $H^0(X;-C-K_X)\cong H^0(X;-L-K_X)\cong H^0(X;-(m-1)H)=0$ and 
$H^1(X;-C-K_X)\cong H^1(X;-L-K_X)\cong H^1(X;L)=0$ (since $|M|$ is  non-special). Therefore, we have
\begin{equation} \label{trisecant}
     H^0(C;{\cal O}_C (C-K_C))\cong H^0(X;-K_X)\cong H^0\left(X;3H-\sum _{i=1}^{3m+3} E_i\right)\cong \C
\end{equation}
and we are done.

Let us prove that $C_m$ is part of a minimal system of generators of $I(X)$. To this end, we first describe $I(X)_2$. Recall that $X$ is the image of the rational map 

$$\varphi :\P^2 \dasharrow \P^N$$
$$(x:y:z)\mapsto (m_0:m_1:\cdots :m_{N-3}:G_x:G_y:G_z),$$

\noindent where $m_i=x^{i_0}y^{i_1}z^{i_2}$ with $i_0+i_1+i_2=m+2$ and $i_0,i_1,i_2\ge 1$. Fix homogeneous coordinates $\omega _0,\ldots ,\omega _N$ in $\P^N$. $I(X)_2$ is generated by
\begin{itemize}
    \item All  $\frac{(m-2)(m-1)m(m+5)}{8}$ quadrics $\omega_i\omega_j -\omega_k\omega_{t}$, $0\le i,j,k,l\le N-3$ satisfying $i_u+j_u=k_u+t _u$ for $u=0,1,2$.
    \item All $\binom{m}{2}$ quadrics $\omega_iH_y-(-1)^{m}\omega_kH_z-\omega_s\omega_{N-2}$ (Analogous quadrics with $\omega_{N-1}$ and $\omega _N$) where $\omega _s=x^{s_0}y^{s_1}z^{s_2}$, $\omega _i=x^{s_0-1}y^{s_1+1}z^{s_2}$ and $\omega _k=x^{s_0-1}y^{s_1}z^{s_2+1}$  with $s_0+s_1+s_2=m+2$, $2\le s_0$ and $1\le s_1,s_2$.
\end{itemize}
 Note that all quadrics in $I(X)_2$ vanish on the plane 
  $\Pi $: $\omega_0=\ldots =\omega_{N-3}=0$ while $C_m$ does not vanish on $\Pi $. Therefore $C_m$ is not a linear combination of the quadrics in $I(X)_2$.
 
 The existence of a trisecant line $\ell $ for $C$ is now clear by the Green-Lazarsfeld criterion mentioned above. If $C=X \cap \P^{N-1}$ then the trisecant line to $C$ is the line $\P^{N-1}\cap \Pi$.

 It remains to see that there is only one cubic in a minimal system of generators of $I(C)$. This follows from the fact that $H^0(C;{\cal O}_C (C-K_C))=\C$ established in \eqref{trisecant}, which means that the span of the unique effective divisor in the bundle ${\cal O}_C (C-K_C)$ is the unique trisecant line. Indeed, assume that $C$ has a unique trisecant line $\ell $ and call $C_1$ the union of $C$ and $\ell $. Then $C_1$ is a canonical singular curve and its ideal is generated by quadrics (Petri’s theorem) which are the quadrics in the ideal of $I(C)$ (In fact, $h^0(C;{\cal O}_C(2)) = h^0(C_1;{\cal O}_{C_1}(2))$ by Riemann-Roch theorem for singular curves (see \cite[Example 18.3.4]{Fu}). Hence, to cut out $C$ from $C_1$, it suffices to consider  one cubic in the ideal of $I(C)$ which is not a linear combination of quadrics.
\endproof
\begin{remark}\label{rem: 5.8} (1) We note that the unexpected curve $C_{R,m}$  in \cite[Theorem 4.6]{Szp19c} can be written
$$\begin{aligned}
C_{R,m}&=xyz\Biggl[(m+1)\Bigl(a(b^m+(-1)^{m+1}c^m)x^{m-1} \\
   &+b(a^m+(-1)^{m+1}c^m)y^{m-1}
   +c(a^m+(-1)^{m+1}b^m)z^{m-1}\Bigr)\\
   &+\sum_{i=1}^{m-2}(-1)^{i+1}\binom{m+1}{i+1}
   \Bigl(a^{m-i}b^{i+1}x^iy^{m-i-1}+ a^{m-i}c^{i+1}x^iz^{m-i-1}+b^{m-i}c^{i+1}y^iz^{m-i-1}\Bigr)\Biggr]\\
   &+a^{m+1}G_x+b^{m+1}G_y+c^{m+1}G_z
\end{aligned}$$

\noindent and this form
does not involve monomials that are a product of $(xyz)^2$ and a monomial of degree $m-4$. 
So, let $J$ be the ideal generated by all coordinates $\omega_i$ corresponding to the monomials in $(xyz)^2\mu_{m-4}$, we can consider the projection $\pi_J$ and consider the map $ \bar\varphi_m=\pi_J \circ\varphi: \P^2\rightarrow \P^{3m-1}$. In this way the ideal of the image $\bar X=\pi_J(X)$ of this map is the elimination ideal of $I(X)$ and so also $I(\bar X)$ is generated by quadrics and a unique cubic which turns out to be the same as in Proposition \ref{prop_image_Fm}. With a simple calculation we get that the number of quadrics decreases to $\frac 1 2(5m^2-m-12)$.
\\
Moreover, $\bar X$ is the result of subsequent projections from external points. So, the degree is preserved and  $deg(\bar X)=deg(X)=m^2+m+1$.
\\
Now we recall that $deg(\bar X)=\sum_i \; h_X(i)$ and, as we know the Hilbert function in degree 1 and 2, we can compute $h_X(1)$ and $h_X(2)$ so we get that $\sum_{i\geq 3}h_X(i)=-(m-2)(m-3)<0$ for any $m\geq 4$. So we deduce that the $h$-vector is not positive and $\bar X$ is not aCM.
Actually, computations with \textit{Macaulay2} supports the
conjecture
that the h-vector is
$$ h=(1, \; 3(m-1),\; 2m^2-7m+9,\;     -\frac 3 2 (m-2)(m-3),\;      \frac 1 2 (m-2)(m-3)).$$
(2) In Proposition 5.7 we have seen that for all $m\ge 4$ the image of the map  $ \varphi: \P^2 \longrightarrow  \P^{N}$ where
  $$\varphi=(xyz\mu_d:G_x:G_y:G_z)$$ is a smooth surface $X$ of degree $m^2+m+1$, its ideal $I(X)$ is generated by quadrics and only one cubic and it is aCM. Using {\em Macaulay2} we have checked what happens with their companion surfaces $X'\subset \P^{3m-1}$. Recall that $X'$ is the image of $\psi: \P^2 \longrightarrow  \P^{3m-1}$ where $$\psi =(a^{m+1},b^{m+1},c^{m+1},a^{m-1}b^2,\ldots,a^2b^{m-1},a^{m-1}c^2,\ldots,a^2c^{m-1},b^{m-1}c^2,\ldots,b^2c^{m-1},$$
  $$a(b^m-(-c)^m,
  b(a^m-(-c)^m),c(b^m-(-a)^m).$$
  Our computation, for $m\le 8$,
  strongly
supports that the following hold: $\psi $ is a birational map onto its image, the ideal of $X'$ is generated by  $\frac 1 2 (5m^2-5m-12)$ quadrics and a unique cubic (unless for $m=4$ where 10 independent cubics are necessary). Moreover, the companion has degree $(m+1)^2$ and with the same arguments as for $\bar X$ we conjecture that the $h$-vector is
  $$ h: (1,\; 3(m-1),\; 2m^2-5m+9,\;     -\frac 3 2 (m-2)(m-3),\;     \frac 1 2 (m-1)(m-6)).$$
Therefore,  $X'$ is not an aCM surface. Indeed, the projective dimension of the coordinate ring of $X'$ is equal to $\codim(X')+1$.

Summarizing,  for $m\leq 8$, we have  that $X\subset \P^N$ is an aCM projection of the Veronese surface $V_{2,m+2}$ in $\P^{\binom{m+4}{2}-1}$ while  its companion surface $X'$ is a non-aCM projection of the Veronese surface $V_{2,m+2}$ in $\P^{3m-1}$. This example makes more intriguing Gr\"obner's problem: To determine when a projection of a Veronese variety $V_{n,d}\subset \P^{\binom{n+d}{d}-1}$ is aCM and when it is not aCM (see \cite{Gr}).
\end{remark}

\begin{remark} If we consider the ideal generated only by the quadrics in $I(\bar X)$ the degree increases by $1$.
Instead  if we consider the ideal generated only by the quadrics in $I(X')$, the degree is still the $(m+1)^2$.

\end{remark}
   Finally, we remark explicitly, that the embedding of $\bar X$ is never too positive.
\begin{corollary}
   The surface $X\subset\P^N$ does not satisfy the $N_2$ condition of Green.
\end{corollary}

\paragraph*{Acknowledgement.}
Our work began during Oberwolfach workshop on Lefschetz Properties in Algebra, Geometry and Combinatorics held in September/October 2020. We thank the MFO for providing, during the COVID-19 crisis, excellent working conditions for collaboration between on site and remote participants of the workshop.

The first two authors have been partially supported by GNSAGA-INDAM.
The third author has been partially supported by PID2020-113674GB-I00. 
The last two authors were partially supported by National Science Centre, Poland, Opus Grant 2019/35/B/ST1/00723. Our collaboration was partially supported by National Science Centre, Poland, Harmonia Grant 2018/30/M/ST1/00148.

We thank the referees for numerous remarks and insights which helped us to improve and clarify our results and their presentation.



\noindent
   Roberta Di Gennaro,\\
   Dipartimento di Matematica e Applicazioni \lq\lq Renato Caccioppoli\rq\rq \\
   Universit\`a degli Studi di Napoli Federico II\\
   Via Cintia, Complesso Universitario Monte Sant'Angelo\\
   80126 Napoli, Italy.
  
\nopagebreak
\noindent
   \textit{E-mail address:} \texttt{digennar@unina.it}\\

\noindent
   Giovanna Ilardi,\\
 Dipartimento di Matematica e Applicazioni \lq\lq Renato Caccioppoli\rq\rq\\
 Universit\`a degli Studi di Napoli Federico II\\
 Via Cintia, Complesso Universitario Monte Sant'Angelo\\
   80126 Napoli, Italy.

\nopagebreak
\noindent
   \textit{E-mail address:} \texttt{giovanna.ilardi@unina.it}\\

\noindent
   Rosa Maria Mir\'o Roig,\\
   Facultat de
  Matem\`atiques i Inform\`atica, 
  Universitat de Barcelona, 
  Gran Via des les Corts Catalanes 585, 
  08007 Barcelona, 
  Spain.
  
\nopagebreak
\noindent
   \textit{E-mail address:} \texttt{miro@ub.edu}\\

\noindent
   Tomasz Szemberg,\\
   Department of Mathematics, Pedagogical University of Cracow,
   Podchor\c a\.zych 2,
   PL-30-084 Krak\'ow, Poland.

\nopagebreak
\noindent
   \textit{E-mail address:} \texttt{tomasz.szemberg@gmail.com}\\

\noindent
   Justyna Szpond,\\
   Institute of Mathematics,
   Polish Academy of Sciences,
   \'Sniadeckich 8,
   PL-00-656 Warszawa, Poland. 
   
\nopagebreak
\noindent
   \textit{E-mail address:} \texttt{szpond@gmail.com}\\

\end{document}